\documentclass[twoside,10pt]{amsart} 
\usepackage{amsfonts,amssymb,latexsym,amsthm, epsfig,amsmath,amscd}
\newcommand{\mf}{\mathfrak}

\newcommand{\lra}{\longrightarrow}

\newcommand{\ol}{\overline}

\newcommand{\ra}{\rightarrow}
\newcommand{\Ra}{\Rightarrow}

\newcommand{\eps}{\varepsilon}

\newcommand{\mbb}{\mathbb}
\newcommand{\tn}{\textnormal}

\newtheorem{de}{Definition}[section]

\newtheorem{pr}[de]{Proposition} 
\newtheorem{tr}[de]{Theorem}
\newtheorem{lm}[de]{Lemma} 

\newtheorem{nt}[de]{Notation} 

\newtheorem{co}[de]{Corollary}

\newcommand{\LMD}{\Lambda}
\newcommand{\lmd}{\lambda}

\def\vp{\rm \vspace{0.2cm}}
\def\M{\rm Max}
\def\m{\mf{m}}

\def\hb{\hfill$\Box$}

\def\Um{\rm Um}
\def\GL{\rm GL}

\def\GQ{\rm GQ}
\def\SL{\rm SL}
\def\EQ{\rm EQ}

\def\E{\rm E}

\def\G{\rm G}
\def\GQ{\rm GQ}

\def\Sp{\rm Sp}

\def\k{\rm K_1}
\def\K{\rm K}
\def\KH{\rm KH}

\def\I{\rm I}

\def\C{\rm C}

\begin{document}
\title{A Note on General Quadratic Groups}
\author{Rabeya Basu\\{\tiny Indian Institute of Science Education and Research, Pune, India}}
\date{}
\maketitle

{\small{\bf Abstract:} We deduce an analogue of 
Quillen--Suslin's local-global principle for the transvection subgroups of the 
general quadratic (Bak's unitary) groups. As an application we revisit the result 
of Bak--Petrov--Tang  on injective stabilization for the ${\k}$-functor of the general 
quadratic groups.} \vp



\section{\large Introduction} 

In this article we are going to discuss three major problems in algebraic ${\K}$-theory 
studied rigorously during 1950's to 1970's, {\it viz.} stabilization problems for the ${\k}$-funtor, 
Quillen--Suslin's local-global principle, and Bak's unitary groups over form rings. 

Let us begin with the stabilization problem for the ${\k}$-functor of the general linear 
groups initiated by Bass--Milnor--Serre during mid 60's. For details, see \cite{BMS}. 
For a commutative ring $R$ with identity, let ${\GL}_n(R)$ be the general linear group, and ${\E}_n(R)$ the subgroup generated by elementary matrices. 
They had studied the following sequence of maps
 $$\frac{{\GL}_n(R)}{{\E}_n(R)}\lra \frac{{\GL}_{n+1}(R)}{{\E}_{n+1}(R)}
\lra \frac{{\GL}_{n+2}(R)}{{\E}_{n+2}(R)}\lra \cdots \cdots = {\k}(R)$$
and posed the problem that when does the natural map stabilize?
That time it was not known that the elementary subgroup is a normal subgroup 
in the general linear group for size $n\ge  3$. 
They showed that the above map is surjective for $n\ge d+1$, 
and injective for $n\ge d+3$, where  $d$ is the Krull dimension of $R$, and 
conjectured that it must be injective for  $n\ge d+2$. In 1969, 
L.N. Vaserstein proved the conjecture 
for an associative ring with identity, where $d$ is the 
stable rank of the ring ({\it cf.} \cite{V}). In 1974, he generalized his result for projective modules. 
After that, in \cite{V2}, he had studied the orthogonal and the 
unitary ${\k}$-functors, and deduced analogue results for those groups, and  showed that for the non-linear cases above mentioned natural map 
is surjective for $n\ge 2d+2$ and 
injectivity for $n\ge 2d+4$. \vp

It was a period when people were working on Serre's problem on 
projective modules, which states that the finitely generated projective 
modules over  polynomial rings over a field are free. If the number 
of variables $n$ in the polynomial ring is 1, then the affirmative answer 
follows from a classical result in commutative algebra. The first non-trivial 
case $n=2$ was proved affirmatively  by C.S. Seshadri in 1958. 
For details see \cite{L}. After 16 years, in 1974, Murthy--Swan--Towber (then Roitman,  
Vaserstein {\it et al.}, for details {\it cf.} \cite{L}) proved the result for algebraically closed fields. The 
general case was proved by Daniel Quillen and Andrei Suslin 
independently in 1976. Suslin proved that for an associative ring $A$, 
which is finite over its center, the elementary subgroup ${\E}_n(A)$ is a normal subgroup of 
${\GL}_n(A)$ for $n\ge 3$. It was first appeared in a paper by M.S. Tulenbaev, for details see \cite{Tu}. 
Quillen, in his proof, introduced a localization technique 
which was one of the main ingredient for the proof of 
Serre's problem (now widely known as Quillen--Suslin Theorem).
Shortly after the original proof, motivated by Quillen's idea, Suslin introduced the 
following matrix theoretic version of the local-global 
principle. We are stating the theorem for a commutative ring with identity, where 
${\M}(R)$ is the maximum spectrum of the ring $R$. The statement is true for 
almost commutative rings, {\it i.e.} rings which are finite over its center. \vp 

{\bf Suslin's Local-Global Principle:} 
{\it Let $R$ be a commutative ring with identity and 
$\alpha(X)\in {\GL}_n(R[X])$ with $\alpha(0)={\I}_n$. If 
$\alpha_{\mf{m}}(X)\in {\E}_n(R_{\mf{m}}[X])$ for every maximal ideal 
$\mf{m}\in {\M}(R)$, then $\alpha(X)\in {\E}_n(R[X])$. } \vp 

In \cite{RB}, we have deduced an analogue of the above statement 
for the transvection subgroups of the full automorphism groups
of global rank at least 1 in the linear case, and at least 
2 for the symplectic and orthogonal cases. In this article we 
generalize the statement for the general quadratic (Bak's unitary) 
groups. The main aim of this article is to show that using this result one can generalize theorems on free 
modules to the classical modules. We shall discuss such technique for the stabilization problems 
for the ${\k}$-functor. In similar manner one can also generalize results for unstable 
${\k}$-groups ({\it cf.} \cite{BBR}, \cite{RB}), 
local-global principle for the commutator subgroups ({\it cf.} \cite{BRK}), 
and may be for the results on congruence subgroup problems, etc. \vp

Let us briefly discuss the historical aspects of the general quadratic 
group defined with the concept of {\it form ring} introduced by Antony Bak 
in his Ph.D. thesis ({\it cf.} \cite{Bak1}) around 1967. For details see \cite{Bak1} and \cite{HV1}. \vp 

We know that a quadratic form on an $R$-module $M$ is a map $q:M\ra R$ such that 
\begin{enumerate}
 \item $q(ax) = a^2q(x)$, \,\,\,\, $a\in R$, \,\,\,\, $x\in M$,
\item  $B_q$ : $M\times M\ra A$ defined by
$B_q(x,y)=q(x+y)-q(x)-q(y)$ is bilinear (symmetric). 
\end{enumerate} 
The map $B_q$ is called the {\it  bilinear form associated to} $q$. If 
$2\in R^{*}$, then for all $x\in M$ we get $q(x)=\frac{1}{2} B_q(x,x)$. 
The pair $(M,q)$ is called the {\it quadratic} $R$-module. 

Suppose $B:V\times V \ra K$ is a bilinear form on a vector space $V$ over a field $K$. We say $B$ is 
\begin{enumerate}
\item  symmetric if $B(u,v)=B(v,u)$, 
\item anti-symmetric if $B(u,v)= - B(v,u)$, and 
\item symplectic (alternating) if $B(u,u)=0$ for all $u\in V$.  
\end{enumerate} 
Now, symplectic $\Ra$ anti-symmetric, and if 
char$(K)\ne 2$, then symplectic $\Leftrightarrow$ anti-symmetric.  
On the other hand, if char$(K)=2$, then anti-symmetric $\Leftrightarrow$ symmetric. 
Also, if char$(K)\ne 2$, then $q$ is a quadratic form if and only if it is 
symmetric bilinear form $B(u,v)=q(u+v)-q(u)-q(v),$ as we get
$\frac{1}{2}B(u,u)=q(u)$. 
The 
study of Dickson--Dieudonne shows that the cases char$(K) = 2$ and char$(K) \ne 2$
differs even at the level of definition. During 1950's--1970's 
it was a problem that how to generalize classical groups by 
constructing theory which does not depend on the invertibility of 2. In 1967,
A. Bak came up with the following: 

First of all, a classical group should be considered as preserving a pair of forms $(B,q)$. 
Secondly, a quadratic form $q$ should take its value not in the ring $R$, but in its 
factor group $R/\Lambda$, where $\Lambda$ is certain additive subgroup of $R$, 
the so called {\it form parameter}. 

In his seminal work ``K-Theory of Forms'', Bak generalized many results for the 
general quadratic groups which were already known for the general linear, symplectic 
and orthogonal groups. In 2003, Bak--Petrov--Tang proved that the injective 
stabilization bound for the ${\k}$-functor of the general quadratic groups over 
form rings is also $2d+4$, like traditional classical groups of even size, ({\it cf.} \cite{Bak2}). 
Later Petrov in his Ph.D Thesis ({\it cf.} \cite{P1}), and Tang (unpublished) independently 
studied the result for the classical modules. But none of those works are published. In this article 
we revisit their results for the general quadratic modules, as an application of the local-global principle 
for the transvection subgroups. But, even though the local-global principle holds for the module finite 
rings, our method is applicable only for the commutative rings. On the other hand, it is possible to apply this 
method to deduce analogue results for any other kind of classical groups, in particular, for the
Bak--Petrov's groups ({\it cf.} \cite{P}), where we 
have analogue local-global principle. 
In this connection, we mention 
that for the structure of unstable ${\k}$-groups of the general quadratic groups, we refer \cite{RB} and \cite{HR}. 
In both the papers it has been shown that the unstable ${\k}$-groups are nilpotent-by-abelian for $n\ge 6$, which generalizes the 
result of Bak for the general linear groups in \cite{Bak}. 
Following our method one can generalize that result for the module case. {\it i.e.} it can be shown that the unstable ${\k}$-groups 
for the (extended) general quadratic modules with rank $\ge 6$ are nilpotent-by-abelian. \vp

\section{\large Preliminaries} 

In this section we shall recall necessary definitions. \vp

\noindent{\bf Form Rings:} 
Let $-:R\ra R$, defined by $a\mapsto \ol{a}$, be an involution on $R$, 
and $\lambda\in C(R)$ = center of $R$ such that  $\lmd \ol{\lmd}=1$. 
We define two additive subgroups of $R$ 
$$\LMD_{\tn{max}}=\{a\in R\,|\, a=-\lmd\ol{a}\} \,\,\,\, \& \,\,\,\,
\LMD_{\tn{min}}=\{a-\lmd\ol{a}\,|\, a\in R\}.$$ 
One checks that $\LMD_{\tn{max}}$ and $\LMD_{\tn{min}}$
are closed under the conjugation operation $a\mapsto \ol{x}ax$ for any
$x\in R$. A $\lmd$-{\it form parameter} on $R$ is an additive subgroup $\LMD$ of $R$
such that $ \LMD_{\tn{min}}\subseteq \LMD\subseteq\LMD_{\tn{max}}$, and 
$\ol{x}\LMD x\subseteq \LMD$ for all $x\in R$. A pair $(R,\LMD)$ is called a {\it form ring}. \vp

\noindent{\bf General Quadratic Groups:} \vp 

Let $V$ be a right $R$-module. By 
${\GL}(V)$ we denote the group of all $R$-linear automorphisms of $V$. 
Throughout the paper we shall consider $V$ as a projective $R$-module. To define 
the {\it general quadratic module} we need following definitions:

\begin{de} \tn{A {\it sesquilinear form} is a map
$f: V\times V \ra R$ such that $f(ua,vb)=\ol{a}f(u,v)b$  for all $u,v\in V$ and $a,b\in R$. }
\end{de} 
\begin{de} \tn{A $\LMD$-{\it quadratic form} on $V$ is a map 
$q : V \ra R/\LMD$ such that $q(v)=f(v,v) + \LMD$.}
\end{de} 

\begin{de} \tn{An {\it associated $\lmd$-Hermitian form} is a map 
$h: V\times V\ra R$ with the property $h(u,v)=f(u,v) + \lmd \ol{f(v,u)}$.} 
\end{de} 
\begin{de} \tn{
A {\it quadratic module} over $(R, \LMD)$  is a triple $(V,h,q)$.}
\end{de}  

The $\lmd$-Hermitian form $h : V\times V \ra R$ induces a map $V \ra$ Hom$_R(V,R)$, given by 
$v\mapsto h(v,-)$. We say that $V$ is {\it non-singular} if $V$ is a projective $R$-module and the Hermitian form 
$h$ is non-singular. 

\begin{de} \tn{
A morphism of quadratic modules over $(R, \LMD)$ is a map 
$\mu : (V,q,\LMD) \ra (V',q',\LMD')$ such that $\mu:V\ra V'$ is $R$-linear, $\mu(\lmd)=\lmd'$, 
and $\mu(\LMD)=\LMD'$.}
\end{de}

\noindent{\bf General Quadratic (Bak's Unitary) Groups:} 
Let $(V, q, h)$ be non-singular quadratic module over $(R, \LMD)$. We define the {\it general quadratic group} as follows:
$${\GQ}(V,q,h) \,\, =\,\,
\{\alpha \in {\GL}(V) \,\,|\,\, h(\alpha u, \alpha v)=h(u,v), \,\, q(\alpha u)=q(v)\}.$$ 
{\it i.e.} the group consisting of all automorphisms which fixes the $\lmd$-Hermitian form and the $\LMD$-quadratic form. 
One observes that the traditional classical groups are the special cases of 
Bak's unitary groups. The central concept is ``form parameter'' due to Bak. 
Earlier version is due to K. McCrimmon which plays an important role in his 
classification theory of Jordan Algebras. He defined for the wider class of 
alternative rings (not just associative rings), but for associative rings it 
is a special case of Bak's concept. For details, see N. Jacobson, Lectures on 
Quadratic Jordan Algebras, TIFR, Bombay 1969, ({\it cf.} \cite{J}) and \cite{HO}. \vp

{\bf Free Case:} 
Let $V$ be a non-singular free right $R$-module of rank $2n$ with ordered basis 
 $e_1,e_2,\ldots, e_n, e_{-n}, \ldots, e_{-2}, e_{-1}$. Consider the 
sesquilinear form $f:V\times V \ra R$ defined by $f(u,v)=\ol{u}_1v_{-1}+\cdots+\ol{u}_nv_{-n}$. 
Let $h$ be the $\lmd$-Hermitian form, and $q$ be the $\LMD$-quadratic form defined by $f$. We get  

$$h(u,v) = \ol{u}_1v_{-1}+\cdots+\ol{u}_nv_{-n} +\lmd \ol{u}_{-n}v_n+\cdots+\lmd\ol{u}_{-1}v_1,$$ 
$$q(u) = \LMD + \ol{u}_1u_{-1}+\cdots+\ol{u}_n u_{-n}.$$
Using the above basis we can identify ${\GQ}(V,h,q)$ with a subgroup of ${\GL}_{2n}(R, \LMD)$ of rank $2n$, say
 ${\GQ}_{2n}(R, \LMD)$. Hence 
$${\GQ}_{2n}(R,\LMD) ~ = ~ \{\sigma\in {\GL}_{2n}(R, \LMD)\,|\, \ol{\sigma}\psi_n
\sigma=\psi_n\},$$
where
$$\psi_n= \begin{pmatrix} 0 & \lmd {\rm I}_n \\ {\rm I}_n &0\end{pmatrix}.$$ 
\noindent{\bf Unitary Transvections $\EQ(V)$:} 
Let  $(V,h,q)$ be a quadratic module over $(R, \LMD)$. 
Let $u,v\in V$ and $a\in R$ be such that $f(u,u)\in \LMD$, $h(u,v)=0$ and $f(v,v)=a$ modulo $\LMD$. 
Then we quote the definition of unitary transvection from (\S 5.1, \cite{B2}). We define
$\sigma=\sigma_{u,a,v}:M\ra M$ defined by 
$$\sigma(x)=x+h(v,x)-v\ol{\lmd}h(u,x)-u\ol{\lmd}ah(u,x).$$

\noindent{\bf Unitary Transvections $\EQ(M)$ in $M = V \perp \mbb{H}(P)$:} 
Let $P$ be a projective $R$-module of rank at least one, and $\mbb{H}(P)$ 
the hyperbolic space. Let $x=(v,p,q)\in M$ for some $v\in V$, $p\in P$, and $q\in P^{*}$. 
For any element $p_0\in P$, $w_0\in V$ and $a_0\in A$ such that $a_0=f(w_0,w_0)$ modulo $\LMD$, 
the above conditions hold, and hence we can define $\sigma_{p_0,a_0,w_0}$ as follows: 
$$\sigma_{p_0,a_0,w_0}(x)=x+p_0h(w_0,x)-w_0\ol{\lmd}h(p_0,x)-p_0\ol{\lmd}a_0h(p_0,x).$$

\noindent{\bf Elementary Quadratic Matrices: (Free Case)} 

Let $\rho$ be the permutation, defined by $\rho(i)=n+i$ for $i=1,\ldots,n$. 
Let $e_i$ denote the column vector with $1$ in the $i$-th position and 0's elsewhere. 
Let $e_{ij}$ be the matrix with $1$ in the $ij$-th position and 0's
elsewhere. For $a\in R$, and  $1\le i, j\le n$, we define 
$$q\eps_{ij}(a)={\rm I}_{2n}+ae_{ij}-\ol{a}e_{\rho(j)\rho(i)} \,\,\tn{ ~~~~for } i\ne j,$$
$$qr_{ij}(a)=\begin{cases} {\rm I}_{2n}+ae_{i\rho(j)}-\lmd\ol{a}e_{j\rho(i)} \,\,\tn{ for }  i\ne j\\
{\rm I}_{2n}+ae_{\rho(i)j}\,\,\,\,\,\,\,\,\,\,\,\ \,\,\,\,\,\,\,\,\,\,\,\,\, \,\,\,\tn{~~ for } 
i=j,\end{cases}$$ 
$$ql_{ij}(a)=\begin{cases}
{\rm I}_{2n}+ae_{\rho(i)j}-\ol{\lmd}\ol{a}e_{\rho(j)i} \,\,\tn{ for } i\ne j\\
{\rm I}_{2n}+ ae_{\rho(i)j} \,\,\,\,\,\,\,\,\,\,\,\ \,\,\,\,\,\,\,\,\,\,\,\,\, \,\,\,\tn{ ~~for } 
i=j.\end{cases}$$
(Note that for the second and third type of elementary matrices, if $i=j$, 
then we get $a=-\lmd\ol{a}$, and hence it forces that $a\in \LMD_{\rm max} (R)$. 
One checks that these
above matrices belong to ${\GQ}_{2n}(R,\LMD)$; {\it cf.}~\cite{Bak1}.) \vp 

\noindent
{\bf n-th Elementary Quadratic Group
${\EQ}_{2n}(R,\LMD)$:} \tn{The subgroup generated by $q\eps_{ij}(a)$, 
$qr_{ij}(a)$ and $ql_{ij}(a)$, for $a\in R$ and $1\le i,j\le n$.} \vp

\noindent{\bf Stabilization map:} There are  standard embeddings:

\begin{center}
${\GQ}_{2n}(R, \LMD) \lra {\GQ}_{2n+2}(R, \LMD)$
\end{center} 
given by 

\begin{center}
$\begin{pmatrix}
  a & b \\
c & d 
 \end{pmatrix} \mapsto 
\begin{pmatrix}
a & 0 & 0 & b \\
0 & 1 & 0 & 0 \\
0 & 0 & 1 & 0 \\
c & 0 & 0 & d 
\end{pmatrix} \tn{\,\,or\,\,} 
\begin{pmatrix}
a & 0 & b & 0 \\
0 & 1 & 0 & 0 \\
c & 0 & d & 0 \\
0 & 0 & 0 & 1 
\end{pmatrix} 
 \tn{\,\,or\,\,} 
\begin{pmatrix}
1 & 0 & 0 & 0 \\
0 & a & b & 0 \\
0 & c & d & 0 \\
0 & 0 & 0 & 1 
\end{pmatrix}.
$\end{center}
Hence we have 
\begin{center}
 ${\GQ}(R,\LMD) = \underset{\lra}\lim\,\, {\GQ}_{2n}(R, \LMD)$.
\end{center} 

It is clear that the stabilization map takes 
generators of ${\EQ}_{2n}(R,\LMD)$ to the generators of 
${\EQ}_{2n+2}(R, \LMD)$. \vp 

\noindent{\bf Commutator Relations:} 
There are standard formulas for the commutators between quadratic 
elementary matrices. For details, we refer \cite{Bak1} (Lemma 3.16). 
In later sections there are repeated use of those relations. \vp

\noindent{\bf Remark:} If $M = R^{2n}$, then under the choice of the above basis 
we get ${\EQ}(M)={\EQ}_{2n}(R,\LMD)$. ({\it cf.} proof of Lemma 2.20 in \cite{BBR}).  \vp 
Hence we have 
\begin{center}
  ${\EQ}(R,\LMD) = \underset{\lra}\lim \,\,{\EQ}_{2n}(R, \LMD)$
\end{center} \vp 
Using analogue of the Whitehead Lemma for the general quadratic groups ({\it cf.} \cite{Bak1}) due to Bak, one gets 
$$[{\GQ}(R,\LMD), {\GQ}(R,\LMD)]=[{\EQ}(R,\LMD), {\EQ}(R,\LMD)]={\EQ}(R,\LMD).$$
Hence 
we define the {\it Whitehead group} of the general quadratic group 
$${\k}{\GQ}={\GQ}(R,\LMD)/{\EQ}(R,\LMD).$$
And, the Whitehead group at the level $m$ 
$${\K}_{1,m}{\GQ}  = {\GQ}_m(R,\LMD)/{\EQ}_m(R,\LMD),$$
where $m=2n$ in the non-linear cases. For classical modules we replace 
$(R,\LMD)$ by $M=V\oplus \mbb{H}(P)$. 

\section{\large Suslin's Local-Global Principle for Tansvection Subgroups} 

In this section we prove analogue of Quillen--Suslin's local-global principle for the transvection subgroups 
of the general quadratic groups. We start with the following splitting lemma. 

\begin{lm} $(${\it cf.} pg. 43-44, Lemma 3.16, \cite{Bak1}$)$ 
\label{key2} 
Let $q_{ij}$ denote any one of the elementary generator of $q\eps_{ij}$, $ql_{ij}$ and $qr_{ij}$ in ${\GQ}_{2n}(R,\LMD)$. 
Then,  for all $x, y \in R$,  $$q_{ij}(x+y)=q_{ij}(x)q_{ij}(y).$$
\end{lm}

We shall need following standard fact.

\begin{lm} \label{key6} 
Let $G$ be a group, and $a_i$, $b_i \in G$, for 
$i = 1, \ldots, n$. Then  for $r_i={\underset{j=1}
{\overset{i}\Pi}} a_j$, we have
${\underset{i=1}{\overset{n}\Pi}}a_i b_i=
{\underset{i=1}{\overset{n}\Pi}}r_ib_ir_i^{-1} 
{\underset{i=1}{\overset{n}\Pi}}a_i.$
\end{lm}
\begin{nt} \tn{ By ${\GQ}_{2n}(R[X], \LMD [X], (X))$ we shall mean the group of all invertible matrices in  ${\GQ}_{2n}(R[X], \LMD [X])$ 
 which are ${\I}_{2n}$ modulo $(X)$. Let $\LMD[X]$ denote the $\lambda$-form parameter on $R[X]$ induced from $(R,\LMD)$, {\it i.e.}, $\lmd$-form parameter 
on $R[X]$ generated by $\LMD$, {\it i.e.}, the smallest form parameter on $R[X]$ containing $\LMD$. 
Let $\LMD_s$ denote the $\lambda$-form parameter on $R_s$ induced from $(R,\LMD)$.}
\end{nt} 
\begin{lm} \label{key1} The group 
${\GQ}_{2n}(R[X],\LMD [X], (X)) \cap {\EQ}_{2n}(R[X], \LMD [X])$
is generated by the elements of the types 
$ \eps q_{ij}(f(X))\eps^{-1}$, where  $\eps \in {\EQ}_{2n}(R, \LMD)$,
$f(X)\in R[X]$ and  $q_{ij}(f(X))$ are
congruent to ${\rm I}_{2n}$ modulo $(X)$.
\end{lm}
{\bf Proof.} 
Let $\alpha(X)\in{\EQ}_{2n}(R[X], \LMD [X])$ be such that 
$\alpha(X)={\rm I}_{2n}$ modulo $(X)$. 
Then we can write $\alpha(X)$ as a product of elements of the 
form $q_{ij}(f(X))$, where $f(X)$ is 
a polynomial in $R[X]$.
We write each $f(X)$ as a sum of a constant term and a 
polynomial which is identity modulo $(X)$. 
Hence by using the splitting property described in Lemma \ref{key2} each 
elementary generator $q_{ij}(f(X))$ 
can be written as a product of two such elementary 
generators with the left one defined on $R$ and the right 
one defined on $R[X]$ which is congruent to ${\rm I}_{2n}$ modulo $(X)$. 

Therefore, we can write $\alpha(X)$ as a product of elementary generators of the form 
$$q_{i j}(f(0))q_{i j}(Xg(X))
 \,\,{\rm for\,\, some }  \,\, g(X)\in R[X]\,\, {\rm with } \,\,g(0)\in R .$$
Now the result follows by using the identity described in Lemma \ref{key6}.  \hb \vp

By repeated application of commutator formulas stated in 
(\cite{Bak1}, pg. 43-44, Lemma 3.16) one gets the following lemma. 

\begin{lm} \label{key3}
Suppose $\vartheta$ is an elementary generator of the general quadratic 
group ${\GQ}_{2n}(R[X], \LMD[X])$. Let $\vartheta$ be 
congruent to identity modulo $(X^{2m})$, for $m>0$. Then, if we conjugate 
$\vartheta$ with an elementary generator of the general quadratic 
group ${\GQ}_{2n}(R, \LMD)$, we get the final matrix 
as a product of elementary generators of the general quadratic  group 
${\GQ}_{2n}(R[X], \LMD[X])$ each of which is 
congruent to identity modulo $(X^m)$. 
\end{lm}

\begin{co} \label{key4} In Lemma \ref{key3} we can take $\vartheta$ as a product of 
elementary generators of the general quadratic group ${\GQ}_{2n}(R[X], \LMD[X])$.
\end{co}

Let us recall following useful and well known fact. We use this fact for the proof of 
{\it dilatation lemma}.
\begin{lm} $(${\it cf.}, \cite{BV}, Lemma 5.1$)$
\label{noeth} 
Let $A$ be Noetherian ring and $s\in A$. Let $s\in A$ and $s\neq 0$. 
Then there exists a natural number 
$k$ such that the homomorphism ${\G}(A,s^kA, s^k\LMD) \ra {\G}(A_s, \LMD_s)$ 
$($induced by the localization homomorphism $A \ra A_s)$ is injective.  
\end{lm}

We recall that one has any module finite ring $R$ as a direct limit of its  finitely generated subrings. Also, 
${\G}(R, \LMD) = \underset{\lra}\lim \, {\G}(R_i, \LMD_i)$, 
where the limit is taken over all finitely generated subring of $R$. 
Thus, one may assume that $C(R)$ is Noetherian. 
 For the rest of this section we shall consider module finite rings $(R, \LMD)$ with identity.\vp

 In \cite{RB}, we have given a proof of the dialation lemma for the general Hermitian 
groups. For the general quadratic modules, the proof is similar, in fact easier. We are giving the sketch of the proof, 
as it is not clearly written any 
available literature. It is almost done in \cite{P1} 
for the odd unitary groups which contains the general quadratic groups. For the hyperbolic unitary groups 
it is done in \cite{BV}.

\begin{lm} \label{Di2}
{\bf (Dilation Lemma: Free Case)} 
Let $\alpha(X)\in {\GQ}_{2n}(R[X], \LMD [X])$, with $\alpha(0)={\rm I}_{2n}$. If 
$\alpha_s(X)\in {\EQ}_{2n}(R_s[X], \LMD_s [X])$ for some non-nilpotent $s\in R$, 
then $\alpha(bX) \in  {\EQ}_{2n}(R[X], \LMD [X])$ 
for $b\in (s^l) {\C}(R)$, $l\gg 0$. 

$($Actually, 
we mean there exists some $\beta(X)\in {\EQ}_{2n}(R[X], \LMD[X])$ such that
$\beta(0)={\rm I}_{2n}$ and $\beta_s(X)=\alpha_s(bX))$.
\end{lm}

{\bf Proof.}  Given that $\alpha_s(X)\in {\EQ}_{2n}(R_s[X], \LMD_s [X])$.  
Since $\alpha(0)={\rm I_{2n}}$, using Lemma \ref{key1}
we can write $\alpha_s(X)$ as a product of the matrices of the form
$ \eps q_{ij}(h(X))\eps^{-1}$, 
where  $\eps \in {\EQ}_{2n}(R_s, \LMD_s)$,
$h(X)\in R_s[X]$ with $q_{ij}(h(X))$ 
congruent to 
${\rm I}_{2n}$ modulo $(X)$. Applying the homomorphism $X\mapsto XT^d$, for $d\gg 0$, from the polynomial ring $R[X]$ to the 
polynomial ring $R[X,T]$, we look on $\alpha(XT^d)$. 
Note that $R_s[X,T]\cong (R_s[X])[T]$. As $C(R)$ is Noetherian, it follows from Lemma \ref{noeth} 
and Corollary \ref{key4} that over the ring 
$(R_s[X])[T]$  we can write 
$\alpha_s(XT^d)$ as a product of elementary generators of the general quadratic group 
such that each of those elementary generator is congruent to 
identity modulo $(T)$. Let $l$ be the maximum of the powers occurring in
the denominators of those elementary generators. Again, as $C(R)$ is Noetherian, 
by applying the homomorphism $T\mapsto s^mT$ for 
$m\ge l$ it follows from Lemma \ref{key4} that over the ring 
$R[X,T]$  we can write 
$\alpha_s(XT^d)$ as a product of elementary generators of general quadratic group 
such that each of those elementary generator is congruent to identity modulo $(T)$.
{\it i.e.} there exists some $\beta(X,T)\in {\EQ}_{2n}(R[X,T], \LMD [X,T])$ such that
$\beta(0,0)={\rm I}_{2n}$ and $\beta_s(X,T)=\alpha_s(bXT^d)$, for some $b\in(s^l)C(R)$.
Finally, the result follows by putting $T=1$.\hb

\begin{lm} \label{Di3}
{\bf (Dilation Lemma: Module Case)} Assume $M=V\oplus \mbb{H}(R)$, where $V$ is a right $R$-module, and 
$\mbb{H}(R)$ is the usual hyperbolic 
space. Let rank of $M$ is $2n$.
Let us denote $M[X]=(V\perp \mbb{H}(R))[X]$. Let $s\in R$ be such that 
$V_s$ is free. Let $\alpha(X)\in \GQ(M[X])$ with 
$\alpha(0)={\rm Id}$. Suppose $\alpha_s(X)\in {\EQ}_{2n}(R_s[X], \LMD_s[X])$. 
Then there exists $\widehat{\alpha}(X)\in {\EQ}(M[X])$ and $l>0$ such that 
$\widehat{\alpha}(X)$ localizes to $\alpha(bX)$ for some $b\in (s^l)$ 
and $\widehat{\alpha}(0)={\rm Id}$. 
\end{lm}
{\bf Proof.} Arguing as in the proof of Proposition 3.1 in \cite{BBR} we can deduce 
the proof. One observes the repeated use of the 
commutator formulas stated in pg. 43, \cite{Bak1}. \hb

\begin{lm} \label{lg2}
{\bf (Local Global Principle for Tansvection Subgroups)} 
Let $M=V\perp\mbb{H}(R)$, where $V$ is as above. 
Let $M[X]=(V\perp\mbb{H}(R))[X]$. 
 Let $\alpha(X)\in {\GQ}(M[X])$ with 
$\alpha(0)={\rm Id}$. Suppose $\alpha_{\mf{m}}(X)\in {\EQ}_{2n}(R_{\mf{m}}[X], \LMD_{\mf{m}}[X])$ 
for every maximal ideal $\mf{m}$ in $R$. 
Then  $\alpha(X)\in {\EQ}(M[X])$.  
\end{lm}
{\bf Proof.} Since $\alpha_{\m}(X)\in {\EQ}_{2n}(R_{\m}[X], \LMD_{\mf m}[X])$ for all
$\m \in {\M}(C(R))$, for each $\m$ there exists $s\in C(R) - \m$ such that 
$\alpha_s(X)\in  {\EQ}_{2n}(R_s[X], \LMD_s[X])$. 
We consider 
a finite cover of $C(R)$, say $s_1+\cdots+ s_r=1$.  Following Suslin's trick, 
let $$\theta(X,T)=
\alpha_s(X+T)\alpha_s(T)^{-1}.$$ 
Then $\theta(X,T)\in {\EQ}_{2n}((R_s[T])[X], (\LMD_s[T])[X])$ and $\theta(0,T)={\I}_{2n}$. 

Since for $l\gg 0$, $\langle s_1^l,\ldots, s_r^l\rangle = R$, we chose 
 $b_1,b_2,\dots,b_r\in C(R)$, with $b_i\in (s^l)C(R),\, l\gg 0$
such that (A) holds and $b_1+\cdots+b_r=1$. 
Then by dilation lemma (Lemma \ref{Di3}), applied with base ring $R[T]$, there exists $\beta^i(X)\in  {\EQ}(M[X], \LMD[X])$ 
(considering $M$ as an $R[T]$ module) such that $\beta^i_{s_i}(X)= \theta(b_iX,T)$. 
Therefore, $\underset{i=1}{\overset{r}\Pi}\beta^i(X)\in {\EQ}(M[X], \LMD[X])$.
 But, $$\alpha_{s_1\cdots s_r}(X)=
\left(\underset{i=1}{\overset{r-1}\Pi} \theta_{s_1\cdots \hat{s_i}\cdots s_r}(b_iX,T)
{\mid}_{T=b_{i+1}X+\cdots +b_rX}\right) \theta_{s_1\cdots \cdots s_{r-1}}(b_rX,0).$$  
Observe that as a consequence of the Lemma \ref{noeth} it follows that 
the map  $${\EQ}(R,s^kR, s^k\LMD) \ra {\E}(R_s, \LMD_s)$$ for $k\in \mbb{N}$, is injective for each $s=s_i$.
As $\alpha(0)={\I}_{2n}$,  we conclude $\alpha(X)\in {\EQ}(M[X])$.  \hb

\section{\large Stabilization of ${\k\GQ}$} 
   
We recall following result of Bak--Petrov--Tang in \cite{Bak2} for free modules of even size. We are stating the theorem for commutative 
ring,s and for this section we shall work for commutative rings with trivial involution. We consider the underline $R$-module 
$M=V\oplus \mbb{H}(R)$, where $V$ is a right $R$-module and $\mbb{H}(R)$ is the hyperbolic space with usual inner product. 
\begin{tr} \label{BPT} $($Bak--Petrov--Tang$)$
Let $(R, \LMD)$ be a commutative form ring with Krull dimention $d$. If $2n\ge {\rm max}(6, 2d+2)$, then 
${\K}_{1,2n}{\GQ}(R, \LMD)$ is a group, the stabilization maps ${\K}_{1,2n-2}{\GQ}(R, \LMD)\lra {\K}_{1,2n}{\GQ}(R, \LMD)$
is surjective, and the maps $${\K}_{1,2n}{\GQ}(R, \LMD)\lra {\K}_{1,2n+2}{\GQ}(R, \LMD)$$ are isomorphism of groups. 
\end{tr} 

As an application of our local-global principle for the transvection subgroup (Theorem \ref{lg2})
we generalize the above stabilization result for the general quadratic modules. Let us first recall the following 
key lemma of Vaserstein ({\it cf.} Corollary 5.4, \cite{V3}). A proof for the absolute case is nicely written in an unpublished 
paper by Maria Saliani ({\it cf.} Theorem 6.1, \cite{MS}). The relative case follows by using the double ring concept, as it is 
done in \cite{G} (Theorem 4.1).

\begin{lm} \label{v3} 
Let $(R, \LMD)$ be an associative ring with identity with Krull dimension $d$, and $I$ be an ideal 
of $R$. Let  $M=V\oplus \mbb{H}(R)$ and $(M, h, q)$ a general quadratic module of rank $2n\ge {\rm max }
(6, 2d+2)$. Then the group of elementary transvection ${\EQ}(M\perp \mbb{H}(R), I)$ acts transitively 
on the set ${\Um}(M\perp \mbb{H}(R), I)$ of unimodular elements 
which are congruent to $(0,\ldots, 1, \ldots, 0)$
modulo $I$. In other words, 
$$\GQ((M\perp \mbb{H}(R)), I)=\EQ ((M\perp \mbb{H}(R)), I)\GQ(V, I).$$
\end{lm}

Also, we have the standard fact.  

\begin{pr} \label{new} 
 Let $M$ and $V$ be as above, and $\Delta\in \GQ(M\perp \mbb{H}(R))$ and $n\ge 2$. If
$\Delta e_{2n}=e_{2n}$, then $\Delta \in \EQ(M\perp \mbb{H}(R))\GQ(M)$. 
\end{pr}
{\bf Proof.}  Proof goes as in Lemma 3.6 in \cite{BR}. \hb  \vp\\
We deduce the following stabilization result for the general quadratic modules: 

\begin{tr} Let $(R, \LMD)$ be an associative ring with identity with Krull dimension $d$. Let 
$M=V\oplus \mbb{H}(R)$ and $(M, h, q)$ a general quadratic module of rank $2n\ge {\rm max }
(6, 2d+2)$. Then, the stabilization map 
 $${\K}_{1,2n}{\GQ}(M)\lra {\K}_{1,2n+2}{\GQ}(M\perp \mbb{H}(R))$$ is isomorphism of groups. 
\end{tr}

{\bf Proof.}
 It follows from the above result of Bak--Petrov--Tang (Theorem \ref{BPT})
that 
for every localization at maximal ideals of $R$ the map is surjective at the level $2d+2$. 
Hence the surjectivity follows by local-global principle, and by the surjectivity result of Bak--Petrov--Tang.

In view of their result, it 
is enough to prove the injectivity for $2n=2d+2$. 
Let $n_1=2n+2$. 
Suppose $\gamma\in \GQ(M)$ is such that $\widetilde{\gamma}=\gamma\perp {\rm Id}$ lies in ${\EQ}(M\perp \mbb{H}(R))$. 
Let $\phi(X)$ be the isotopy between $\widetilde{\gamma}$ and identity. 
Viewing $\phi(X)$ as a matrix it follows that $v(X)=\phi(X) e_{n_1}$ is in 
${\Um}(M\perp \mbb{H}(R))[X])$, and $v(X)=e_{n_1}$ modulo $(X^2-X)$. 
It follows from Lemma \ref{v3} that over $R[X]$ we get $\sigma(X)\in {\EQ}((M\perp \mbb{H}(R))[X])$ such that 
$\sigma(X)^tv(X)=e_{n_1}$ and $\sigma(X)={\rm Id}$ modulo $(X^2-X)$. Therefore, 
$\sigma(X)^t\phi(X) e_{n_1}=e_{n_1}$. Then by  Lemma \ref{new} over $R[X]$ we can write 
$\sigma(X)^t\phi(X)=\psi(X)\widetilde{\phi}(X)$ for some $\psi(X)\in {\EQ}((M\perp \mbb{H}(R))[X])$ and 
$\widetilde{\phi}(X)\in \GQ(M[X])$. 

Since $\sigma(X)={\rm Id}$ modulo $(X^2-X)$, $\widetilde{\phi}(X)$ is an isotopy  between $\gamma$ 
and identity. Therefore,  after localization at a
maximal ideal $\mf{m}$, the image $\widetilde{\phi}_{\mf{m}}(X)$ is stably elementary 
for every maximal ideal $\mf{m}$ in ${\M} (R)$. Hence by the stability theorem of Bak--Petrov--Tang 
(Theorem \ref{BPT}) for the free modules, it follows that 
$\widetilde{\phi}_{\mf{m}}(X)\in {\EQ}(2d+2, R_{\mf{m}}[X])$. Since $\phi(0)={\rm Id}$, we get
$\widetilde{\phi}(0)={\rm Id}$. Hence by the above local-global principle for the transvection subgroups 
(Theorem \ref{lg2}) it follows that
 $\widetilde{\phi}(X)\in {\EQ}(M[X])$. Hence $\gamma = \widetilde{\phi}(1) \in {\EQ}(M)$. 
This proves that the map at the level $2d+4$ is injective. \vp

{\bf Acknowledgment:} My sincere thanks to Sergey Sinchuk for many useful discussions. I am highly 
grateful to Nikolai Vavilov and Alexei Stepanov for their kind invitation  to visit Chebyshev Laboratory, 
St. Petersburg, Russia, where I could finish my long pending work on this manuscript. I am thankful 
to St. Petersburg State University and IISER Pune for supporting my visit.

{\tiny
 
\addcontentsline{toc}{chapter}{Bibliography} 
{\scshape Indian Institute of Science Education and Research (IISER Pune)},\\
{\scshape Dr. Homi Bhabha Road, Pashan,}\\
{\scshape Pune, Maharashtra 411008, India.} \\
{\tiny Email: rabeya.basu@gmail.com, rbasu@iiserpune.ac.in}}


\begin{thebibliography}{99} 
\bibitem{Bak1} {\scshape A. Bak}; ${\K}$-Theory of Forms. Annals of Mathematics Studies, 
98. Princeton University Press, Princeton, N.J. University of Tokyo Press, 
Tokyo (1981). 
\bibitem{Bak} {\scshape  A. Bak}; Nonabelian ${\K}$-theory: the nilpotent class of 
${\k}$ and general stability. ${\K}$-Theory  4  (1991),  no. 4, 363--397. 
\bibitem{BV} {\scshape A. Bak, N. Vavilov}; 
Structure of hyperbolic unitary groups. I. Elementary subgroups. Algebra Colloq. 7 (2000), no. 2, 159--196. 
\bibitem{Bak2} {\scshape A. Bak, V. Petrov, G. Tang}; Stability for Quadratic ${\k}$. 
${\K}$-Theory 29 (2003), 1--11. 
\bibitem{BBR} {\scshape A. Bak, R. Basu, R.A. Rao}; 
Local-Global Principle for Transvection Groups. Proceedings of The American Mathematical Society 138 (2010), no. 4, 
1191--1204.
\bibitem{BMS} {\scshape H. Bass, J. Milnor, J.-P. Serre}; 
Solution of the congruence subgroup problem for $\SL_n$ $(n\ge 3)$ and $\Sp_{2n}$ $(n\ge 2)$. 
Inst. Hautes $\acute{\E}$tudes Sci. Publ. Math. No. 33 (1967) 59--137.
\bibitem{B} {\scshape H. Bass}; Algebraic K-Theory. Benjamin, New York-Amsterdam (1968). 
\bibitem{B2} {\scshape H. Bass}; Unitary algebraic ${\K}$-theory.  Algebraic K-theory, III: 
Hermitian ${\K}$-theory and geometric applications (Proc. Conf., Battelle 
Memorial Inst., Seattle, Wash., 1972). Lecture Notes in Mathematics, 
Vol. 343, Springer, Berlin (1973), 57--265.
\bibitem{B1} {\scshape H. Bass}; Quadratic modules over polynomial rings, Contributions to Algebra
 (a collection of papers dedicated to Ellis Kolchin), Academic Press, New York (1977).
\bibitem{BRK} {\scshape R. Basu, R.A. Rao, Reema Khanna};  On Quillen's local global principle. Commutative algebra and algebraic geometry, 17–30, Contemp. Math., 390, Amer. Math. Soc., Providence, RI, 2005.
\bibitem{BR} {\scshape R. Basu, R.A. Rao};  
Injective stability for ${\k}$ of the orthogonal group.  Journal of Algebra  323  (2010),  no. 2, 393--396.
\bibitem{RB} {\scshape R. Basu};  
Local-global principle for general quadratic and general Hermitian groups and the nilpotence of 
${\KH}_1$. Zap. Nauchn. Sem. S.-Peterburg. Otdel. Mat. Inst. Steklov. (POMI) 452 (2016), Voprosy Teorii Predstavleniĭ Algebr i Grupp. 30, 5–31.
\bibitem{G} {\scshape A. Gupta};  Optimal injective stability for the symplectic K1Sp group. 
J. Pure Appl. Algebra 219 (2015), no. 5, 1336--1348.
\bibitem{HR} {\scshape R. Hazrat}; Dimension theory and nonstable ${\k}$ of quadratic modules.  ${\K}$-Theory  27  
(2002),  no. 4, 293--328. 
\bibitem{HV1} {\scshape R. Hazrat, N. Vavilov}; Bak's work on ${\K}$-theory of rings. On the occasion of his 65th 
birthday. Journal of ${\K}$-Theory 4 (2009), no. 1, 1--65.
\bibitem{HO} {\scshape A.J. Hahn, O.T. O'Meara}; The Classical groups and ${\K}$-theory. With a foreword by 
J. Dieudonné. Grundlehren der Mathematischen Wissenschaften [Fundamental Principles of Mathematical Sciences], 
291. Springer-Verlag, Berlin, (1989).
\bibitem{J} {\scshape N. Jacobson}; Lectures on Quadratic Jordan algebras, Tata Istitute of Fundamental Research, 
Bombay, (1969). 
\bibitem{KOP}  {\scshape V.I. Kopeiko}; The stabilization of Symplectic groups 
over a polynomial ring. Math. USSR. Sbornik 34 (1978), 655--669.
\bibitem{L} T.Y. Lam; Serre's problem on projective modules. 
Springer Monographs in Mathematics. Springer-Verlag, Berlin, 2006.
%\bibitem{MVV} B.A. Magurn, W. Van der Kallen, L.N. Vaserstein; Absolute stable rank and Witt cancellation 
%for noncommutative rings. Invet. math. 91 (1988), 525--542.  
\bibitem{P} {\scshape V.A. Petrov}; Odd unitary groups. 
(Russian) Zap. Nauchn. 
Sem. S.-Peterburg. Otdel. Mat. Inst. Steklov. (POMI) 305 (2003), Vopr. Teor. Predst. Algebr. i Grupp. 10, 195--225, 
241; translation in J. Math. Sci. (N. Y.) 130, no. 3 (2005), 4752--4766. 
\bibitem{P1} {\scshape V.A. Petrov}; Overgroups of classical groups. Ph.D. thesis, Saint--Petersburg State University (2005), pp 129. 
\bibitem{MS} {\scshape Maria Saliani}; On the stability of the unitary group. Available on the web; Homepage for 
Prof. M.-A. Knus. 
\bibitem{SK} {\scshape A.A. Suslin, V.I. Kopeiko}; Quadratic modules and the orthogonal group over polynomial rings. (Russian) Modules and representations. Zap. Naučn. Sem. 
Leningrad. Otdel. Mat. Inst. Steklov. (LOMI) 71 (1977), 216–250, 287.
\bibitem{Tu} {\scshape M.S. Tulenbaev}; Schur multiplier of a group of elementary matrices of finite order, 
Zapiski Nauchnykh Seminarov Leningradskogo Otdeleniya Matematicheskogo Instuta im V.A. Steklova Akad. Nauk  
SSSR, Vol. 86 (1979), 162--169
\bibitem{V} {\scshape L.N. Vaserstein}; On the Stabilization of the general Linear group 
over a ring.  Mat. Sbornik (N.S.) 79 (121) 405--424 (Russian); English 
translated in Math. USSR-Sbornik. 8 (1969), 383--400. 
\bibitem{V2} {\scshape L.N. Vaserstein}; Stabilization of Unitary and Orthogonal Groups 
over a Ring with Involution. Mat. Sbornik, Tom 81 (123) (1970) no. 3, 307--326.
\bibitem{V3} {\scshape L.N. Vaserstein}; Stabilization for Classical groups over rings. 
(Russian)  Mat. Sb. (N.S.)  93 (135)  (1974), 268--295, 327. 
\bibitem{vas2} {\scshape L.N. Vaserstein}; On the normal subgroups of ${\rm GL}\sb{n}$ over a ring. 
Algebraic $K$-theory, Evanston (1980) (Proc. Conf., Northwestern Univercisy, 
Evanston, Ill., (1980)), pp. 456--465, Lecture Notes in Mathematics, 854, Springer, 
Berlin-New York, (1981). 
\end{thebibliography}
\end{document}